\documentclass[12pt]{article}
\usepackage[usenames,svgnames,dvipsnames]{xcolor}
\usepackage{adjustbox}
\usepackage{amsmath}
\usepackage{amssymb}
\usepackage{xcolor}
\usepackage{graphicx}
\usepackage{overpic}
\usepackage[normalem]{ulem}
\usepackage{ arydshln}
\usepackage[square,numbers]{natbib}
\usepackage{mathtools}
\usepackage{hyperref}
\usepackage{doi}
\usepackage{booktabs}
\usepackage{arydshln}
\usepackage{subfigure}
\usepackage{rotating}
\usepackage{pdflscape}
\usepackage{wasysym}
\usepackage{multirow}
\usepackage{cancel}
\usepackage{array}
\usepackage{comment}

\newcommand{\sspace}[1]{S_{#1}}
\newcommand{\mtrx}[1]{\boldsymbol{\mathrm{#1}}}
\newcommand{\vctr}[1]{\boldsymbol{#1}}
\newcommand{\trans}[1]{#1^T}
\newcommand{\mXi}{\mtrx{\Xi}}
\newcommand{\mZeta}{\mtrx{Z}}

\newcommand{\vxi}{{\vctr{\xi}}}
\newcommand{\vi}{\vctr{i}}
\newcommand{\vj}{{\vctr{j}}}
\newcommand{\vk}{\vctr{k}}

\newcommand{\cc}{\text{cc}}
\newcommand{\fcc}{\text{fcc}}
\newcommand{\bcc}{\text{bcc}}
\newcommand{\hex}{\text{h}}
\newcommand{\NDccII}{\mXi_{\cc 2}}
\newcommand{\NDccIII}{\mXi_{\cc 3}}
\newcommand{\NDfcc}{\mXi_{\fcc}}
\newcommand{\NDbcc}{\mXi_{\bcc}}
\newcommand{\NDqc}{\mXi_{\text{qc}}}
\newcommand{\NDiii}{\mXi_{3}}

\newcommand{\tabref}[1]{Table~\ref{#1}}
\newcommand{\secref}[1]{Section~\ref{#1}}

\definecolor{blue}{rgb}{0,0.3,.7}

\newcommand{\minho}[1]{\textcolor{orange}{#1}}

\newcommand{\vol}[1]{\operatorname{vol}\left(#1\right)}
\newcommand{\R}{\mathbb{R}}

\newcommand{\Z}{\mathbb{Z}}

\newcommand{\bxs}{{box spline}}
\newcommand{\Bxs}{{Box spline}}
\newcommand{\BxS}{{Box Spline}}

\newcommand{\jorg}[1]{\textcolor{purple}{#1}}

\newcommand{\Mxzp}{\Mcc{21}}

\newcommand{\MhexII}{\Mhex{20}}

\newcommand{\matid}{\mtrx{I}}
\newcommand{\Zfcc}{\Z_{\fcc}}
\newcommand{\Zbcc}{\Z_{\bcc}}
\newcommand{\Zhex}{\Z_{\hex}}
\newcommand{\ZAn}[1]{\mathbb{A}_{#1}}
\newcommand{\ZAndual}[1]{\mathbb{A}_{#1}^*}
\newcommand{\ZDn}[1]{\mathbb{D}_{#1}}
\newcommand{\ZDndual}[1]{\mathbb{D}_{#1}^*}

\newcommand{\Gfcc}{\mGen_{\fcc}}
\newcommand{\Gbcc}{\mGen_{\bcc}}
\newcommand{\Ghex}{\mGen_{\hex}}
\newcommand{\GAn}[1]{\mtrx{A}_{#1}}
\newcommand{\GAndual}[1]{\mtrx{A}_{#1}^*}
\newcommand{\GDn}[1]{\mtrx{D}_{#1}}
\newcommand{\GDndual}[1]{\mtrx{D}_{#1}^*}
\newcommand{\perm}{\pi}
\newcommand{\mGen}{\mtrx{G}}
\newcommand{\ZG}{\Z_{\mGen}}
\newcommand{\Half}{\frac{1}{2}}

\newcommand{\Mcc}[1]{M_{\text{c}#1}}
\newcommand{\Mhex}[1]{M_{\text{h}#1}}
\newcommand{\Mfcc}[1]{M_{\text{f}#1}}
\newcommand{\Mbcc}[1]{M_{\text{b}#1}}

\newcommand{\mP}{\mtrx{P}}

\newcommand{\mL}{\mtrx{L}}

\newcommand{\opspan}{\operatorname{span}}
\newcommand{\rank}{\operatorname{rank}}

\newcommand{\figref}[1]{Fig.~\ref{#1}}

\newcommand{\SymGrp}[1]{\mathcal{SG}\left(#1\right)}
\newcommand{\PntSet}[2]{\mathcal{PS}\left(#1,#2\right)}
\newcommand{\DirSet}[2]{\mathcal{DS}\left(#1,#2\right)}

\def\makehsquare#1#2#3{
    \dimen0=#1\advance\dimen0 by -#3
    \vrule height#1 width#2 depth0pt \kern-#2
    \vrule height#1 width#1 depth-\dimen0 \kern-#1
    \vrule height#2 width#1 depth0pt \kern-#3
    \vrule height#1 width#3 depth0pt
}

\def\boxo{\mathop{\makehsquare{6pt}{1.2pt}{.05pt}\kern.05em}}
 
\def\subboxo{\mathop{\makehsquare{4pt}{1pt}{.05pt}}}

\begin{document}
\title{A Practical \BxS\ Compendium}

\author{Minho Kim, J{\"o}rg Peters}

\date{\today}
\maketitle

\begin{abstract}
\Bxs s provide smooth spline spaces as shifts of a single generating function on a lattice and so generalize tensor-product splines. Their elegant theory is laid out 
in classical papers and a summarizing book.
This compendium aims to succinctly but exhaustively survey symmetric low-degree \bxs s
with special focus on two and three variables.
Tables contrast the lattices, supports, 
analytic and reconstruction properties,
and list available implementations and code.
%
\end{abstract}

\def\wid{\0.2\linewidth}
\begin{figure}[ht]
\begin{center}
\begin{tabular}{cc}
\subfigure[shifts of the univariate `hat' function
on $\Z$]{
\label{fig:slab-1d}
\includegraphics{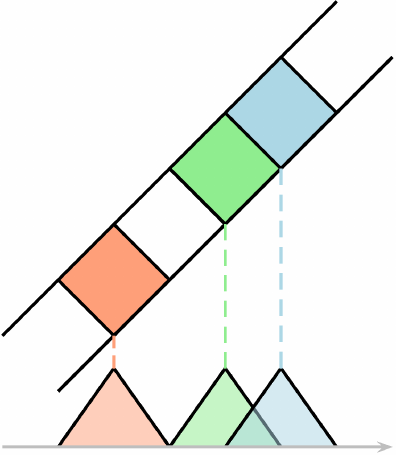}
}
&
\subfigure[shifts of 
 the bivariate `hat' function 
$M_{h10}$ on $\Zhex$]{
\label{fig:slab-2d}
\includegraphics{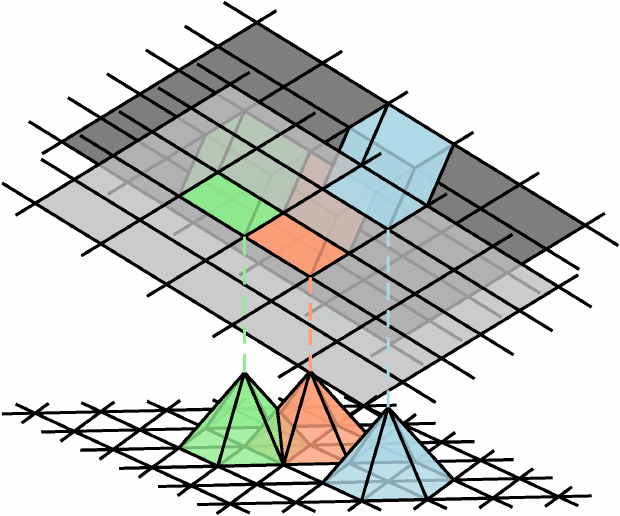}
}
\end{tabular}
\end{center}
\caption[]{
\Bxs s 
as a projection
of  $n$-dimensional boxes \cite{Kim2010Symmetric}.
}
\label{fig:projDef}
\end{figure}

\section{Introduction}
\label{sec:intro}

As a generalization of uniform polynomial tensor-product splines, and with the beautiful interpretation as a projection
of a higher-dimentional box partition 
\cite{Sommerfeld1904,deBoor1983Approximation,Prautzsch2002Box,Kim2010Symmetric},
see \figref{fig:projDef},  \bxs s
have repeatedly commanded the attention of researchers seeking an elegant foundation for differentiable function spaces on low-dimensional lattices. 
Notably, \bxs s provide the regular prototypes for generalized uniform polynomial
subdivision algorithms \cite{catmull1978recursively,DeRose1998Subdivision,peters2008subdivision} and have been advocated for reconstructing signals on
non-Cartesian lattices,
see \secref{sec:reconstruction}.
This compendium summarizes the 
latest findings for \bxs\ spaces with emphasis on $d=2$ and $d=3$ variables
and
 \emph{symmetric}
\bxs s, i.e.\ \bxs s that have
at least the symmetry of their domain lattice.
The aim is to provide a succinct overview, via tables and illustrations,  of the properties, 
literature and
 computational tools
 and code, and 
to characterize each \bxs's  efficiency
in terms of smoothness, polynomial reproduction, support size and polynomial degree.

\section{Lattices and \bxs s}
\label{sec:Def}
We refer to \citet{Conway2013Sphere}
for a general treatment of lattices and their symmetry groups, beyond the needs of the compendium.
\paragraph{Lattices and Direction Sets}

Given the integer grid $\Z^d$,
any non-singular
$d\times d$ \emph{generator matrix} $\mGen$
defines a lattice $\ZG:=\mGen\Z^d$.
The \emph{symmetry group} $\SymGrp{\ZG}$ of $\ZG$,
represented as an orthogonal matrix group,
consists of all
orthogonal transformations
that leave $\ZG$ invariant:
\[
    \SymGrp{\ZG}:=
    \left\{
        \mL\in\R^{d\times d}:\trans{\mL}\mL=\matid_d\text{ and } \forall \vj\in\ZG\  \mL\vj\in\ZG
    \right\}.
\]
where $\matid_d$ is the $d\times d$ identity matrix. 

\begin{table}[!htb]
    \caption{Five domain lattices for $d=2,3$.
    }
    \label{tab:domain.lattices}
    \centering
    \begin{tabular}{ccccc}
        \toprule
         dim. & name & symbol & generator matrix & 
         $\#\SymGrp{*}$ \\
         \midrule
         \multirow{2}{*}{$2$} 
         & Cartesian & $\Z^2$ & $\matid_2$ & $8$ \\
         \cmidrule{2-5}
            & hexagonal & $\Zhex$ & $    \Ghex :=
            \frac{1}{2}
            \begin{bsmallmatrix*}[r]
                1 & 1 \\
                -\sqrt{3} & \sqrt{3}
            \end{bsmallmatrix*}
            $
    & 
    $12$
   \\ 
        \cmidrule{1-5}
        \multirow{6}{*}{$3$} 
        & Cartesian & $\Z^3$ & $\matid_3$ & $48$ \\
        \cmidrule{2-5}
        & FCC & \multirow{2}{*}{$\Zfcc$} & \multirow{2}{*}{$\Gfcc:=
    \begin{bsmallmatrix*}[r]
        0   &   1   &   1   \\
        1   &   0   &   1   \\
        1   &   1   &   0   \\
    \end{bsmallmatrix*}$} & \multirow{2}{*}{$48$} \\
        & (face-centered cubic) \\
        \cmidrule{2-5}
        & BCC & \multirow{2}{*}{$\Zbcc$} & \multirow{2}{*}{$ \Gbcc:=
    \begin{bsmallmatrix*}[r]
        -1   &   1   &   1   \\
        1   &   -1   &   1   \\
        1   &   1   &   -1   \\
    \end{bsmallmatrix*}$} & \multirow{2}{*}{$48$} \\
    & (body-centered cubic) \\
         \bottomrule\multicolumn{5}{l}{\scriptsize  
         $\#\mathcal{S}$ is the cardinality of the set 
    $\mathcal{S}$, $\matid_d$ the $d\times d$ identity matrix.}
    \end{tabular}
\end{table}
\begin{figure}[h]
    \begin{center}
        \subfigure[$\Z^2$]{
        \label{fig:Z2.dir}
     \includegraphics[width=.3\textwidth]{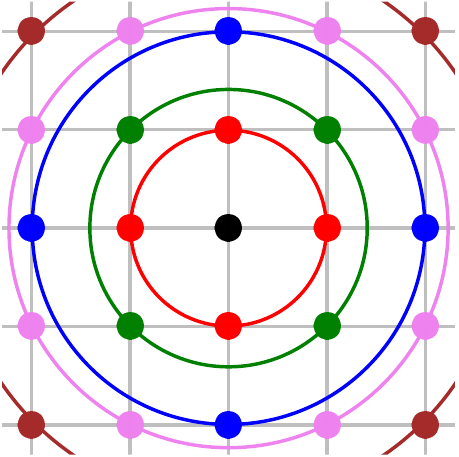}
        }
        \hskip1cm
        \subfigure[$\Zhex$]{
        \label{fig:hex.dir}
    \includegraphics[width=.3\textwidth]{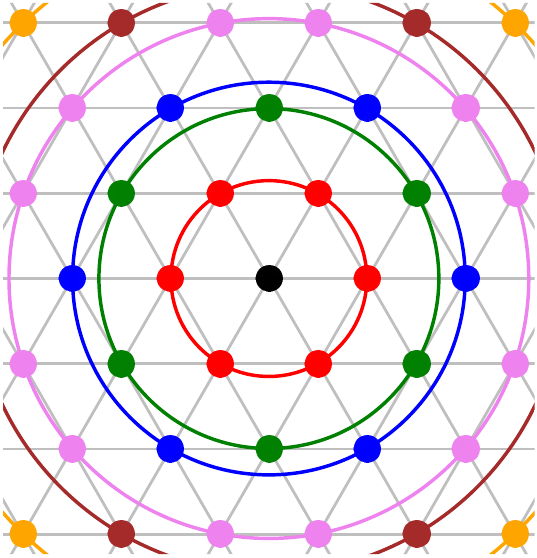}
        }
    \end{center}
    \caption{
Stratifying 2D lattice points 
by distance to the origin {\Large$\bullet$}
such that each shell corresponds to a  direction set
        $\DirSet{\ZG}{k}$, 
        \textcolor{Red}{$k=1$},
   \textcolor{Green}{$k=2$},
   \textcolor{Blue}{$k=3$},
   \textcolor{Magenta}{$k=4$}, $\ldots$ 
    }
    \label{fig:dir.2d}
%
    \centering
    \subfigure[$\Z^3$]{
    \label{fig:cc.dir}
    \includegraphics[width=.3\textwidth]{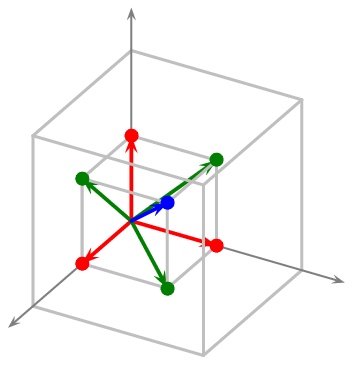}
    }
    \subfigure[$\Zfcc$]{
    \label{fig:fcc.dir}
    \includegraphics[width=.3\textwidth]{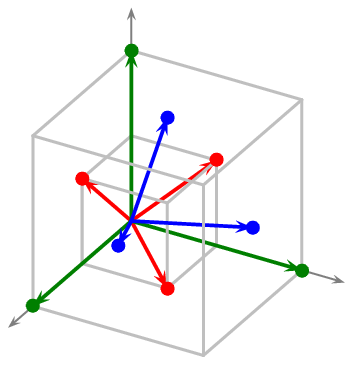}
    }
    \subfigure[$\Zbcc$]{
    \label{fig:bcc.dir}
    \includegraphics[width=.3\textwidth]{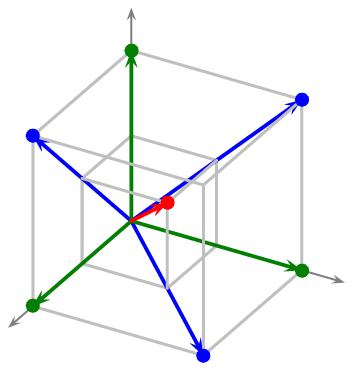}
    }
    \caption{
    Stratifying 3D 
    direction vectors corresponding to
        direction sets
        $\DirSet{\ZG}{k}$, 
   \textcolor{Red}{$k=1$},
   \textcolor{Green}{$k=2$},
   \textcolor{Blue}{$k=3$}.
 \tabref{tab:lattice.directions} lists 
    coordinates.
    }
    \label{fig:dir.3d}
\end{figure}
In the plane (2D)  and 3-space (3D), five  lattices are known
for their high symmetries.
They are listed in \tabref{tab:domain.lattices}.
To enumerate \bxs s, we collect the lattice direction
vectors $\vj \in \ZG$ 
into \emph{direction sets}
$\DirSet{\ZG}{k}$
consisting of one vector and its images under the symmetry group of the lattice. The index $k$ 
is assigned by non-decreasing vector
length, 
see \figref{fig:dir.2d} and \figref{fig:dir.3d},
which is unique for $k\le 3$, the cases of interest. (For $k>3$, 
multiple direction sets can lie in the same 
spherical shell \cite{Conway2013Sphere},
e.g.\ $(5,0)$ and $(4,3)$ in $\Z^2$.)
Since $-\vj=\mGen(-\vi)$ and $-\vi\in\Z^d$
if $\vi\in\Z^d$,
for each $\vj=\mGen\vi\in\ZG$ also $-\vj\in\ZG$,
we list only one of $\vj$ and $-\vj$
in $\DirSet{\ZG}{k}$.

\begin{figure}[ht]
\begin{center}
\subfigure[
$
\mXi := 
\begin{bsmallmatrix}
    1 & 0 \\
    0 & 1
\end{bsmallmatrix}
$
]{
\label{fig:zp.step.1}
\includegraphics[width=.3\textwidth]{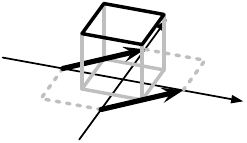}
}
\subfigure[$
\mXi := 
\begin{bsmallmatrix}
    1 & 0 & 1\\
    0 & 1 & 1
\end{bsmallmatrix}
$]{
\label{fig:zp.step.2}
\includegraphics[width=.3\textwidth]{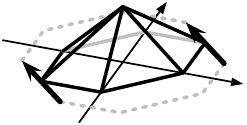}
}
\subfigure[$
\mXi := 
\begin{bsmallmatrix}
    1 & 0 & 1 & -1\\
    0 & 1 & 1 & 1
\end{bsmallmatrix}
$]{
\label{fig:zp.step.3}
\includegraphics[width=.3\textwidth]{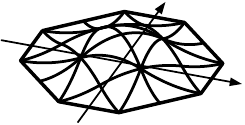}
}
\end{center}
\caption[]{
\Bxs s 
via convolution in the directions
(columns) of $\mXi$
on $\Z^2$.
}
\label{fig:convolution}
\end{figure}

\paragraph{\BxS s}
Given a domain lattice $\ZG$,
direction vectors $\vxi\in\ZG$
can 
be collected into a $d\times m$ 
\emph{direction matrix} 
$\mXi$ to 
 define the centered \bxs\ $M_{\mXi}$ recursively,
 starting with the characteristic function
 $\chi_{\mXi\subboxo^d}$ on the (half-open) parallelepiped  $\mXi\boxo^d$, $\boxo:=\left[-\Half,\Half\right)$, see  \cite{Chui1988Multivariate,deBoor1993Box} 
 and
\figref{fig:convolution}:
\begin{equation}
\label{eq:def.conv}
    M_{\mXi}:=
    \begin{cases}
        \displaystyle 
        \int_{-\Half}^{\Half} M_{\mXi\backslash\vxi}
            \left(\cdot-t\vxi\right){\rm d}t & \text{if }d<m,\ \vxi\in\mXi,  \\
        \displaystyle 
        \frac{|\det{\mGen}|}{|\det{\mXi}|}\chi_{\mXi\subboxo^d}
                &\text{if }d=m \text{ and }\det\mXi\ne 0.
    \end{cases}
\end{equation}
\def\wid{\0.2\linewidth}

The centered \bxs\ is invariant under exchange of columns
or multiplication of a column by -1:
$    M_{\mXi_1}=M_{\mXi_2}
$
     if and only if there exists 
a `signed permutation' matrix $\mP$
that can permute and/or change sign of a coordinate,
 such that
$\mXi_1=\mXi_2\mP$.
Moreover, 
since for any linear map $\mL$, see \cite[page 11]{deBoor1993Box},
\begin{equation}
\label{eq:lin.map}
    M_{\mXi}=
    |\det{\mL}|M_{\mL\mXi}(\mL\cdot),
\end{equation}
many properties for centered \bxs s on the Cartesian lattice  $\Z^d$ transfer directly to 
$\ZG$
by a linear change of variables $\mGen$.


Let $\mXi\in \mGen\Z^{d\times m}$
with $\rank{\mXi}= d$, $M_{\mXi}$
the corresponding \bxs, and 
$\sspace{\mXi} := \opspan(M_{\mXi}(\cdot-\vj))$ the space
of its shifts over the lattice.
Then $M_{\mXi}$ and $\sspace{\mXi}$ 
have the following properties:
\begin{enumerate}
\item
 $M_{\mXi}$ is
non-negative
and its
 shifts over $\ZG$ sum to 1:
due to  the  factor $|\det\mGen|$ in
 \eqref{eq:def.conv}
\[
    \sum_{\vj\in\ZG}M_{\mXi}(\cdot-\vj)=1.
\]
\item
 The support of $M_{\mXi}$  is 
 $\mXi\boxo^d$, i.e.\ the set sum of the 
 vectors in $\mXi$.
\item
 $M_{\mXi}$ is  piecewise polynomial of total degree $m-d$.
    \item 
$M_{\mXi} \in C^{r-2}$. i.e.\ $r-2$ times continuously differentiable, 
where $r$ is the minimal number of columns that need to be removed from  $\mXi$ to obtain a matrix whose columns do not span $\R^d$.
\item
$\sspace{\mXi}$ reproduces all polynomials of degree  $r-1$.
\item
The $L^p$ approximation order
of $\sspace{\mXi}$ is $r$
\cite[page 61]{deBoor1993Box}, i.e.\ for all sufficiently smooth $f$ there exists a sequence $c:\ZG\mapsto \R$
such that
:
\begin{equation}
\label{eq:a.o}
\left\| f - \sum_{\vj\in\ZG} c(\vj) M_{\mXi}((\cdot-\vj)/h)\right\|_p
= O(h^r),  \quad h<1. 
\end{equation}
\item
$\sspace{\mXi}$ forms a basis (the shifts are linearly independent) if and only if 
all square nonsingular submatrices of $\mXi$ are unimodular, i.e.,
$|\det\mZeta|=1$ for all $\mZeta\subset \mXi$ where $\mZeta\in\R^{d\times d}$ 
\cite[{page 41}]{deBoor1993Box}.
\item
With $\vol{\mXi\boxo^d}$ denoting the volume of the support of $M_{\mXi}$,  
the number of coefficients on $\ZG$
required to evaluate a spline value
is $\vol{\mXi\boxo^d}/|\det{\mGen}|$, \cite[{page 36}]{deBoor1993Box}.

\end{enumerate}
The \emph{symmetry group} of $M_{\mXi}$ is defined
analogous to the symmetry group of a lattice:
\[
    \SymGrp{M_{\mXi}}:=
    \left\{
        \mL\in\R^{d\times d}:\trans{\mL}\mL=\matid_d\text{ and }
        M_{\mXi}=M_{\mXi}(\mL\cdot)
    \right\}.
\]
A centered \bxs\ $M_{\mXi}$ on the domain lattice $\ZG$ is 
\emph{symmetric} if it has the same or more symmetries
than $\ZG$:
$\SymGrp{\ZG}\subset \SymGrp{M_{\mXi}}$.
(The centered \bxs\ defined by
$\mXi := 
\begin{bsmallmatrix}
    1 & 1 \\
    0 & 1
\end{bsmallmatrix}
$ 
is not symmetric:  its
symmetry group is $\{\matid_2,-\matid_2\}$,
but the symmetry group of $\Z^2$ has the cardinality $8$
of the signed permutation group.)
If $\vxi \in \DirSet{\ZG}{k}$ is a column of $\mXi$
then 
all directions of $\DirSet{\ZG}{k}$ must be columns in $\mXi$ to make $M_{\mXi}$ symmetric.
This can be seen as follows.
For any $\vxi\in\ZG$, let
$
        \mXi  :=
        \left\{
            \mL\vxi:\mL\in\SymGrp{\ZG}
        \right\}.
$
Then for any $\mL \in\SymGrp{\ZG}$,
the set of directions 
$\mXi$
 equals
the set 
    $\mL\mXi$
    and $|\det{\mL}|=1$
so that by \eqref{eq:lin.map}
$
    M_{\mXi}=
    |\det{\mL}|M_{\mL\mXi}(\mL\cdot)
    =
        M_{\mXi}(\mL\cdot).
$
That is,
    $M_{\mXi}$ is symmetric.
 It suffices to include either $\vxi$ or $-\vxi$ into $\mXi$
    since for any $\vxi\in\DirSet{\ZG}{k}$
    \[
        \int_{-1/2}^{1/2}f(\cdot-t\vxi)dt
        =
        \int_{-1/2}^{1/2}f(\cdot-t(-\vxi))dt
        =
        \int_0^{1/2}f(\cdot-t\vxi)dt
        +\int_0^{1/2}f(\cdot-t(-\vxi))dt.
    \]


\section{Choice of direction vectors}
\label{sec:choice}
The  algebraic and differential geometric properties of \secref{sec:Def}
imply that the
efficiency of a \bxs\ space is closely related to the choice of direction vectors in the construction of the \bxs\ and 
favors the 
vectors to be
\begin{itemize}
\item[-]
snapped to a grid:
this guarantees that the approximation order can be maximal.
(In the extreme case,
the shifts of $M_{[1/2]}$ on $\Z$
do not sum to 1. The shifts 
of $M_{[1,1/2]}$ on $\Z$ 
form a partition of 1, but 
a spline in $\sspace{[1,1/2]}$
has intervals where the spline is constant and cannot match  linear functions.)
\item[-]
 short:  since longer vectors result in larger support 
 and
 more vectors are required to achieve symmetry,
 increasing the degree.
\item[-]
uniformly distributed: for the same degree,
uniformity increases the continuity and approximation order.
(For example, see \tabref{tab:bivariate}, the bi-linear B-spline $\Mcc{20}$ and the ZP element $\Mcc{11}$ have degree $2$, 
but 
both the continuity and the approximation order of $\Mcc{11}$ is higher by one than those of  $\Mcc{20}$.)
\item[-]
in $\DirSet{\ZG}{1}$:
for the five lattices, direction sets with $k>1$ yield $\mXi$ that are not unimodular, and so 
the \bxs\ shifts are not linearly independent
 \cite{deBoor1993Box}.
\end{itemize}
Uniform distribution on a lattice is in 
competition with shortness since equi-distribution of directions requires inclusion of farther lattice points.
\begin{table}[!htb]
\caption{
   The direction sets
    of the five domain lattices in \tabref{tab:domain.lattices}:
   repeating directions are grayed out.
  Numbers in the parentheses denote the cardinality of 
corresponding direction set, cf.\ \figref{fig:dir.2d} and \figref{fig:dir.3d}.
    }
    \label{tab:lattice.directions}
    \centering
    \begin{tabular}{ccrcrcrcr}
        \toprule
       \multirow{2}{*}{lattice} &
       \multicolumn{8}{c}{
       $\DirSet{\ZG}{k}$ 
       }
       \\
       &
       \multicolumn{2}{c}{\textcolor{Red}{$k=1$}} &
   \multicolumn{2}{c}{\textcolor{Green}{$k=2$}} &
   \multicolumn{2}{c}{\textcolor{Blue}{$k=3$}} &
   \multicolumn{2}{c}{\textcolor{Magenta}{$k=4$}}
   \\
        \midrule
         $\Z^2$ & $\NDccII$ & ($2$) & $\NDqc$ & ($2$) & 
         \textcolor{gray}{$2\NDccII$} & ($2$)
    & $\{\perm(2,\pm1)\}$ & ($4$)
        \\
         $\Zhex$ & $\Ghex\NDiii$   & ($3$) & 
                $\Ghex  
                \begin{bsmallmatrix*}[r]
                    2 & -1 & -1 \\
                    1 & 1 & -2
                \end{bsmallmatrix*}
                $ & ($3$) & \textcolor{gray}{$\Ghex(2\NDiii)$} & ($3$)
                &
                $\Ghex
                    \begin{bsmallmatrix*}[r]
                        1 & 2 & 3 & 3 & 2 & 1 \\
                        3 & 3 & 2 & 1 & -1 & -2
                    \end{bsmallmatrix*}
                $
                & ($6$)
                \\
         $\Z^3$ & $\NDccIII$ & ($3$) &  $\NDfcc$ & ($6$) & $\NDbcc$ & ($4$) & \textcolor{gray}{$2\NDccIII$} & ($3$) \\
         $\Zfcc$ & $\NDfcc$ & ($6$) & $2\NDccIII$ & ($3$) & $\{\perm(2,\pm1,\pm1)\}$ & ($12$) & \textcolor{gray}{$2\NDfcc$ } & ($6$) \\
         $\Zbcc$ & $\NDbcc$ & ($4$) & $2\NDccIII$ & ($3$) & $2\NDfcc$ & ($6$) & $\{\perm(3,\pm1,\pm1)\}$  & ($12$) \\
         \bottomrule
         \multicolumn{9}{l}{\scriptsize{
          $\{\perm{(x_1,x_2,\ldots,x_d)}\}$
    is the set 
    of vectors
    generated by permuting the coordinates $x_i$.}
    }\\
    \multicolumn{9}{l}{\scriptsize{
    E.g.\ $\{\perm(2,\pm1)\} = \{(2,1),(2,-1),(1,2),(-1,2)\}$}
    }
    \end{tabular}
    
\end{table}

\tabref{tab:lattice.directions}
lists the direction sets
for the bivariate and trivariate 
domain lattices 
of \tabref{tab:domain.lattices}
in terms of the matrices
(see \figref{fig:dir.2d} and \ref{fig:dir.3d}):
\begin{align*}
 d=2:\quad  & 
    \NDccII := \matid_2,\qquad
     \NDqc :=
     \begin{bmatrix*}[r]
        1 & -1  \\
        1 & 1 
    \end{bmatrix*},
    \qquad
    \NDiii :=
     \begin{bmatrix*}[r]
      1 &0 &-1 \\
      0 &1 & -1
    \end{bmatrix*},
   \\
    d=3:\quad  & 
    \NDccIII:=\matid_3,\ 
    \NDfcc:=
    \begin{bmatrix*}[r]
    1 & -1 & 1 & 1 & 0 & 0 \\
    1 & 1 & 0 & 0 & 1 & -1 \\
    0 & 0 & 1 & -1 & 1 & 1
    \end{bmatrix*},
    \quad
    \NDbcc:=
    \begin{bmatrix*}[r]
    -1 & 1 & 1 & -1  \\
    1 & -1 & 1 & -1  \\
    1 & 1 & -1 & -1
    \end{bmatrix*},
\end{align*}
where the subscripts are to remind of
Cartesian (cc2, cc3)
quincunx (qc), 
3 directions,
FCC, and BCC directions, respectively.


\begin{figure}[ht]
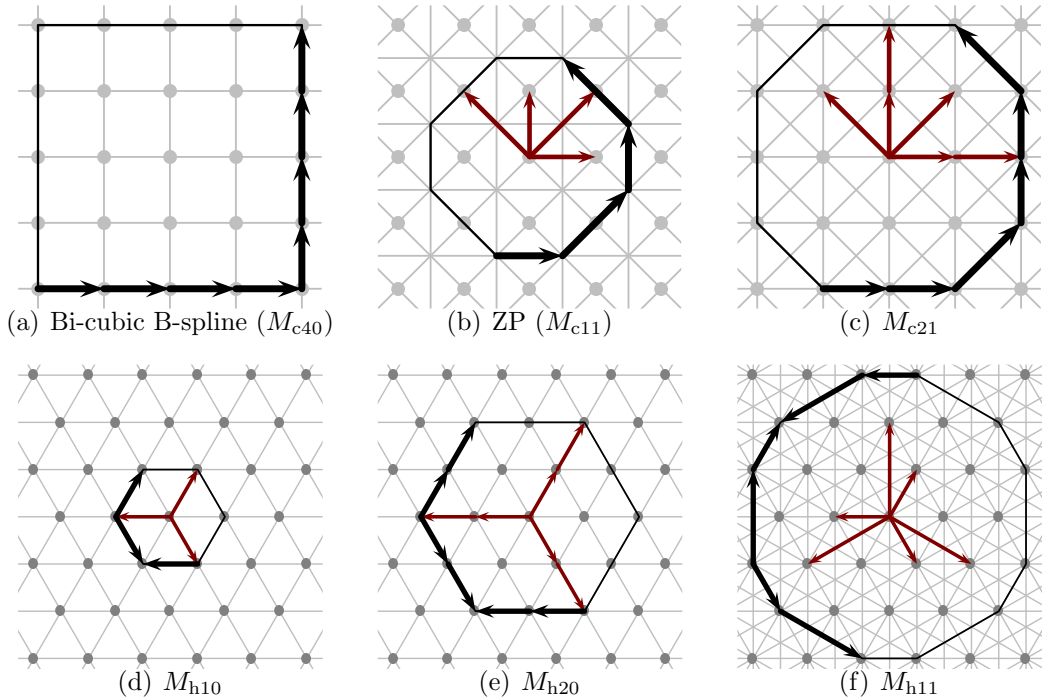

    \centering
    \begin{tabular}{ccc}
        \subfigure[Bi-cubic B-spline ($\Mcc{40}$)]{
			\label{fig:biBspline}	
        \adjustimage{height=.26\textwidth,margin=0pt 0pt,valign=M}{Figures/Z2/Mc40-dir-supp}	
		} 
    &
         \subfigure[ZP ($\Mcc{11}$)]{
			\label{fig:Mc11.dir.supp}	
        \adjustimage{height=.26\textwidth,margin=0pt 0pt,valign=M}{Figures/Z2/Mc11-dir-supp}	
		} 
  &
        \subfigure[ $\Mcc{21}$]{
			\label{fig:Mc21.dir.supp}	
       \adjustimage{height=.26\textwidth,margin=0pt 0pt,valign=M}{Figures/Z2/Mc21-dir-supp}
		} 
    \\
    \subfigure[$\Mhex{10}$]{
			\label{fig:Mh10.dir.supp}
       \adjustimage{height=.26\textwidth,margin=0pt 0pt,valign=M}{Figures/hex/Mh10-dir-supp}	
		} 
  &
         \subfigure[$\Mhex{20}$]{
			\label{fig:Mh20.dir.supp}	
                \adjustimage{height=.26\textwidth,margin=0pt 0pt,valign=M}{Figures/hex/Mh20-dir-supp}
		} 
  &
        \subfigure[$\Mhex{11}$]{
			\label{fig:Mh11.dir.supp}	
       	   \adjustimage{height=.26\textwidth,margin=0pt 0pt,valign=M}{Figures/hex/Mh11-dir-supp}
		} 
    \end{tabular}
    \caption{
    Directions (arrows) and supports 
    (polygons with black edges)
    of
    select bivariate \bxs s
    with polynomial pieces delineated by
    knot lines (gray lines).}
    \label{fig:biDirSup}
\end{figure}
\begin{table}[!htb]
\caption{
Bivariate symmetric \bxs s up to degree $6$.
$\Mcc{n0}$ is the 
tensor-product B-spline, 
$\Mcc{11}$ is the Zwart-Powell element, 
$\Mcc{21}$ is the extended 
6-direction ZP element,
and
$\Mhex{10}$
the hat function. 
The continuity is $C^{r-2}$ with $r$ defined by 4. of \secref{sec:Def}.
}
\label{tab:bivariate}
\begin{center}
\begin{tabular}{cccccccccc}
\toprule
\multirow{2}{*}{lattice} & 
\multicolumn{2}{c}{direction sets } & 
\multirow{2}{*}{degree} & differentiability &
stencil & 
\multirow{2}{*}{reference}
\\
& 
1 & 2  & & $r{-}2 =$ & size 
\\
\midrule
\multirow{5}{*}{$\Z^2$}
&  $n$ & $0$  & $2n{-}2$ & $n{-}2$ & $n^2$ 
& \cite{deBoor1978Practical} 
\\
& $1$ & $1$ & $2$ & $1$ & $7$ 
&\cite{Zwart1973Multivariate,Powell1974Piecewise,peters1997Simplest}
\\
 &  $2$ & $1$  & $4$ &  $2$ & $14$ &  
\cite{Chui1991Algorithms,Peters2004Combining,Lai2007Spline}
 \\
 & $3$ & $1$ & $6$ & $3$ & $23$
 \\
 & $2$ & $2$ & $6$ & $4$ & $28$
\\
\cmidrule{1-7}
\multirow{2}{*}{$\Zhex$}
& $n$ & $0$  & $3n{-}2$ & $2n{-}2$ & $3n^2$  
& 
\cite{Frederickson1970Triangular,Frederickson1971Generalized,Loop1987Smooth,Lai2007Spline}
\\
& $1$ & $1$ & $4$ & $3$ & $24$ 
&  \\
\bottomrule
\end{tabular}
\end{center}
\end{table}
\section{Bivariate \bxs s}
\label{sec:bi}
Since the 
third direction set in \tabref{tab:lattice.directions}
of $\Z^2$ and $\Zhex$ already repeat
 the first,
we restrict the list of bivariate \bxs s in 
\tabref{tab:bivariate} 
to $\DirSet{\ZG}{k}$ for $k<3$, as illustrated in 
\figref{fig:biDirSup}.
We could skip $k=3$ and consider 
the \bxs\ defined by 
$\cup_{k=1,2,4}\DirSet{\Z^2}{k}$
with $2+2+0+4=8$ directions,
but the corresponding \bxs\ 
has a large support and  degree  $8-2=6$, while the resulting $C^5$ continuity is unlikely to match any generic application needs. 
Similarly, 
the \bxs\ defined by 
$\DirSet{\Z^2}{4}$
yields a \bxs\ of degree 2 but of support size 24, whereas the ZP spline 
$M_{c11}$ has the same smoothness for support size 7.
\begin{figure}[ht]
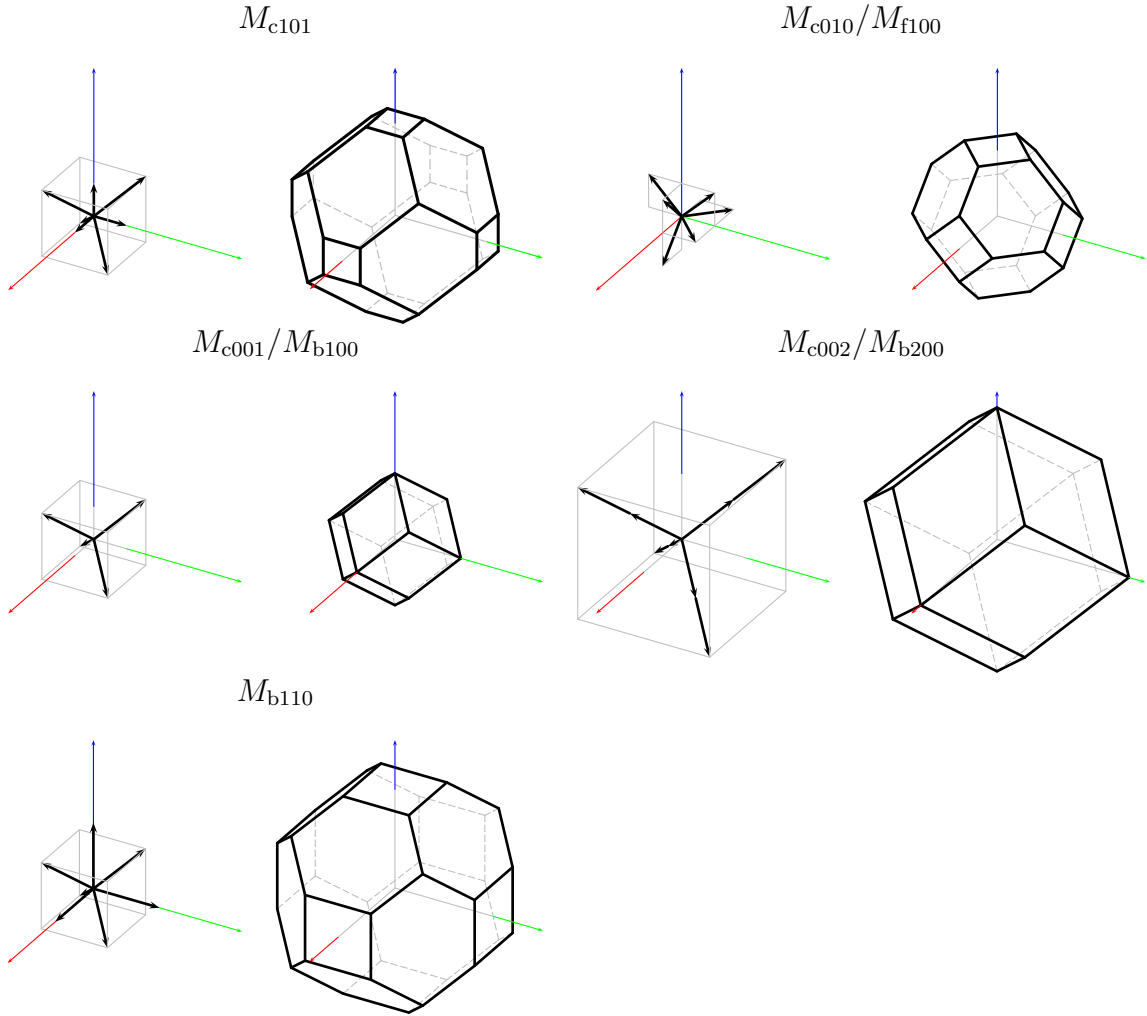

\begin{center}
\begin{tabular}{rrrrrrrr}
    \multicolumn{2}{c}{$\Mcc{101}$} 
    & \multicolumn{2}{c}{$\Mcc{010}$/$\Mfcc{100}$}
    \\
    \adjustimage{width=.2\textwidth,valign=T}{Figures/cc7/cc7-dir} 
    & \adjustimage{width=.216\textwidth,valign=T}{Figures/cc7/cc7-supp}
    & \adjustimage{width=.2\textwidth,valign=T}{Figures/fcc6/fcc6-dir}
    & \adjustimage{width=.2\textwidth,valign=T}{Figures/fcc6/fcc6-supp}
    \\
    \multicolumn{2}{c}{$\Mcc{001}$/$\Mbcc{100}$} 
&    \multicolumn{2}{c}{$\Mcc{002}$/$\Mbcc{200}$}   
    \\
    \adjustimage{width=.2\textwidth,valign=T}{Figures/bcc4/bcc4-dir} 
    & \adjustimage{width=.2\textwidth,valign=T}{Figures/bcc4/bcc4-supp}
    & \adjustimage{width=.216\textwidth,valign=T}{Figures/bcc8/bcc8-dir-repeated}
    & \adjustimage{width=.24\textwidth,valign=T}{Figures/bcc8/bcc8-supp}
    \\
    \multicolumn{2}{c}{$\Mbcc{110}$} 
    \\
    \adjustimage{width=.2\textwidth,valign=T}{Figures/bcc7/bcc7-dir}
    & \adjustimage{width=.228\textwidth,valign=T}{Figures/bcc7/bcc7-supp}
\end{tabular}
\end{center}
\label{fig:3dsup}
\caption{Directions and supports
of select 
trivariate \bxs s.
}
\end{figure}

\begin{table}[!htb]
\def\wid{0.15\linewidth}
\caption{
Trivariate symmetric \bxs s up to degree $9$.
}
\label{tab:trivariate}
\begin{center}
\begin{tabular}{cccccccc}
\toprule
\multirow{2}{*}{lattice} 
& \multicolumn{3}{c}{direction sets} & \multirow{2}{*}{degree} & 
differentiability &
stencil &
note  \\
 & 1 & 2 & 3 & & 
 $r{-}2=$
& size & / reference\\
\midrule
\multirow{9}{*}{$\Z^3$}
&$n$ & $0$ & $0$ & $3n{-}3$ & $n{-}2$ & $n^3$ 
&
B-splines
\cite{deBoor1973Spline}
\\
& $1$ & $1$ & $0$ & $6$ & $3$ & $87$ 
& \cite{Entezari2001Optimal}
\\
& $2$ & $1$ & $0$ & $9$ & $4$ & $172$
\\
& $1$ & $0$ & $1$ &$4$ & $2$ & $53$ 
& \cite{Peters1994C2,Peters1997Box,Entezari2006Extensions,Kim2021Fast}
\\
& {$1$} & {$0$} & {$2$} & {$8$} & {$4$} & {$249$}
\\
& $2$ & $0$ & $1$ & $7$ & $4$ & $120$
\\
& {$0$} & {$n$} & {$0$} & {$6n{-}3$} & {$3n{-}2$} & {$32n^3$}
&\cite{Kim2023Volume}\\
& $0$ & $1$ & $1$ & $7$ & $5$ & $216$
\\
& $0$ & $0$ & $n$ & $4n{-}3$ & $2n{-}2$ & $16 n^3$ & 
\\
\cmidrule{1-8}
\multirow{5}{*}{$\Zfcc$}
& {$n$} & {$0$} & {$0$} & {$6n{-}3$} & {$3n{-}2$} & {$16n^3$} 
& 
\cite{Kim2008Box,Kim2013GPU}
 
\\
& $1$ & $1$ & $0$ & $6$ & $3$ & $86$   
& \cite{Entezari2001Optimal}$^\dagger$
\\
& $1$ & $2$ & $0$ & $9$ & $4$ & $228$
\\
&   $0$ & $n$ & $0$ & $3n{-}3$ & $n{-}2$ & $4n^3$ & B-splines
\\
& $0$ & $0$ & $1$ & $9$ & $7$ & $784$ & 
\\
\cmidrule{1-8}

\multirow{9}{*}{$\Zbcc$}
& $n$ & $0$ & $0$ & $4n{-}3$ & $2n{-}2$ & $4n^3$ 
& \cite{Entezari2004Linear,Kim2010Symmetric}
\\
& {$2$} & {$1$} & {$0$} & {$8$} & {$4$} & {$106$}
\\
& $1$ & $1$ & $0$ & $4$ & $2$ & $30$  
&  \cite{Kim2013Quartic}
\\
& $1$ & $2$ & $0$ & $7$ & $4$ & $92$
\\
& {$1$} & {$0$} & {$1$} & {$7$} & {$5$} & {$200$}
\\
&  $0$ & $n$ & $0$ & $3n{-}3$ & $n{-}2$ & $2n^3$ &  
 B-splines \cite{Csebfalvi2006Prefiltered}
\\
& {$0$} & {$1$} & {$1$} & {$6$} & {$3$} & {$174$}
\\
& {$0$} & {$2$} & {$1$} & {$9$} & {$4$} & {$344$}
\\
& {$0$} & {$0$} & {$n$} & {$6n{-}3$} & {$3n{-}2$} & {$64n^3$} \\
\bottomrule\multicolumn{8}{l}{\scriptsize $^\dag$
The 
\bxs\ proposed in \cite{Entezari2001Optimal} is 
a sibling of $\Mfcc{110}$
built from
the direction matrix 
$\begin{bmatrix}\NDfcc & \NDccIII \end{bmatrix}$.
Since $\NDccIII$ 
}
\\
\multicolumn{8}{l}{\scriptsize
do not snap to $\Zfcc$, the
resulting approximation order is lower
than $\Mfcc{110}$}.

\\
\end{tabular}
\end{center}
\end{table}
Denoting by $n_k$ the number of
repetitions of the $k$th direction set,
 the \bxs\ 
 \begin{equation*}
 \text{on } \Z^2 \text{ are named 
 } \Mcc{n_1n_2} \text{ and those on
  } \Zhex\ \text{   are named } \Mhex{n_1n_2}.
\end{equation*}
%
\tabref{tab:bivariate} leaves out 
direction sets of the form $(0,n)$ and $(1,n)$ for $\Z^2$, 
since their properties do not improve on $(n,0)$ and $(n,1)$, respectively and result in a larger support.
Analogously, $(0,n)$ is omitted for $\Zhex$.
We note that the options for 
$C^1$ continuity are 
$M_{c30}$ (9), $M_{c11}$ (7),
with the stencil sizes
 listed in parentheses.
For 
$C^2$ continuity the options are
$M_{c40}$ (16), $M_{c21}$ (14), and
$M_{h20}$ (12).
%
The only linearly independent
symmetric \bxs s are $\Mcc{n0}$,
i.e.\
the B-splines on
$\Z^2$, and
$\Mhex{n0}$ on 
 $\Zhex$.
(Other  linearly independent \bxs s, such as the three-direction \bxs\ on 
 $\Z^2$
\cite{deBoor1983Bivariate}, are not symmetric.
)
The stencil size explains why several \bxs s have not been investigated in detail.

\section{Trivariate \bxs s}

Analogous to the bivariate case,
denoting by $n_k$ the number of
repetitions of the $k$th direction set,
the \bxs s on $\Z^3$, $\Zfcc$, and $\Zbcc$ are named 
\begin{equation*}
 \Mcc{n_1n_2n_3}, \Mfcc{n_1n_2n_3}, 
 \text{ and } \Mbcc{n_1n_2n_3}
\end{equation*}
in \tabref{tab:trivariate}. 
Fourth 
 direction vectors are not used since, e.g.\ for  $\Mbcc{*}$,
 they
are typically too long and too many.
The only symmetric linearly-independent \bxs s are $\Mcc{n00}$, the B-splines on 
 $\Z^3$,
$\Mfcc{n00}$ on 
 $\Zfcc$, and 
$\Mbcc{n00}$
 on 
 $\Zbcc$.
(There are additional linearly independent asymmetric \bxs s like 
 four-direction \bxs s on $\Z^3$).
 That is $M_{*n00}$ are the only \bxs s that form a basis.
 Listing the support sizes in parentheses, 
 the 
 $C^1$ \bxs s are
$M_{c300}$ (27), $M_{c010}$ (32),
$M_{f100}$ (16), $M_{f030}$ (108),
$M_{b030}$ (54), $M_{b001}$ (64)
and the 
$C^2$  \bxs s are
$M_{c400}$ (64), $M_{c101}$ (53),
$M_{c002}$ (128),
$M_{f040}$ (256), 
$M_{b200}$ (32), $M_{b110}$ (30),
$M_{b040}$ (128).
Due to their small support and the degree listed in
square brackets, $M_{c300}$ [6], $M_{c010}$ [3],
$M_{f100}$ [3] (see \figref{fig:3dsup}) stand out as efficient for $C^1$ and 
$M_{b200}$ [5], $M_{b110}$ [4]  for $C^2$.


\section{Multi-variate \bxs s}

\begin{table}[!htb]
    \caption{
    The first and second direction sets of 
    the five main 
    lattices
    $\ZAn{d}:=\GAn{d}\Z^d$,
    $\ZAndual{d}:=\GAndual{d}\Z^d$,
    $\ZDn{d}:=\GDn{d}\Z^d$, and
    $\ZDndual{d}:=\GDndual{d}\Z^d$.
   For ease of notation,
    opposite directions
    $
    \{\pm\vj:\vj\in\DirSet{\ZG}{k}\}
    $ 
    are enumerated  
    and the directions of $\ZAndual{d}$ are scaled by $(d+1)$  and 
   those
    of
    $\ZDndual{d}$ by $2$.
    As in \cite{Conway2013Sphere}, $a^\alpha$ abbreviates $\alpha$-fold  repeating 
    entries $a,\dots,a$.
    }
    \label{tab:DS.nvariate}
    \centering
    \begin{tabular}{cccc}
        \toprule
         lattice & dim. & $k=1$ & $k=2$ \\
         \midrule
         $\Z^d$ & $d{\ge}2$ & $\{\perm(\pm1,0^{d-1})\}$ & $\{\perm((\pm1)^2,0^{d-2})\}$ \\
         \cmidrule{1-4}
         \multirow{2}{*}{$\ZAn{d}$} & $d{=}2$ & $\{\perm(1,-1,0)\}$ & $\{\pm\perm(2,-1,-1)\}$ \\
         &  $d{>}2$ & $\{\perm(1,-1,0^{d-1})\}$ & $\{\perm(1^2,(-1)^2,0^{d-4})\}$
         \\
         \cmidrule{1-4}
         $\ZAndual{d}$ & $d{\ge}2$
         & $\{\pm\perm(d,(-1)^d)\}$
         & $\{\pm\perm((d-1)^2,(-2)^{d-1})\}$
         \\
         \cmidrule{1-4}
            \multirow{2}{*}{$\ZDn{d}$} & $d{=}3$
         & $\{\perm(\pm1,\pm1,0)\}$
         & $\{\perm(\pm2,0,0)\}$
         \\
          & $d{>}3$
         & $\{\perm((\pm1)^2,0^{d-2})\}$
        & $\{\perm(\pm2,0^{d-1})\}\cup\{\perm((\pm1)^4,0^{d-4})\}$
         \\
         \cmidrule{1-4}
            \multirow{5}{*}{$\ZDndual{d}$}
         & $d{=}2,3$ 
            & $\{((\pm1)^d)\}$
            & $\{\perm(\pm2,0^{d-1})\}$
            \\
          &  $d{=}4$
         & $\{((\pm1)^4)\}\cup\{\perm(\pm2,0^3)\}$
        & $\{\perm((\pm2)^2,0^2)\}$
        \\
       & $4{<}d{<}8$
 & $
\{\perm(\pm2,0^{d-1})\}
$ 
& 
$ \{((\pm1)^d)\} $
\\
         & $d{=}8$
        &
$ \{\perm(\pm2,0^{d-1})\} $
& 
$ 
\{((\pm1)^d)\}\cup\{\perm((\pm2)^2,0^{d-2})\}
$ 
\\
 & $d{>}8$
 & 
 $ \{\perm(\pm2,0^{d-1})\} $
 & 
$
\{\perm((\pm2)^2,0^{d-2})\}
$
\\
\bottomrule
 \end{tabular}
\end{table}

The five lattices in two and three variables are
instances of $d$-dimensional lattices, $d>3$
whose detailed definition can be found in
\cite{Kim2010Symmetric,Kim2011Symmetric}.
The generator matrices of the four lattices other than $\Z^d$ \cite{Conway2013Sphere}
are as follows.
\[
    \GAn{d}:=
        \begin{bmatrix*}[r]
            -1 &  &    & \\
            1 & -1 &    & \\
             & 1 & \ddots &  &  \\
            &  & \ddots & -1 \\
             &   &    & 1 & -1 \\
             &   &    &  & 1 \\
        \end{bmatrix*},
    \GAndual{d} :=
        \frac{1}{d+1}
        \begin{bmatrix*}[r]
            d & -1 & \cdots & -1 & -1 \\
            -1 & d & \cdots & -1 & -1 \\
            \vdots & \vdots & \ddots & \vdots & \vdots \\
            -1 & -1 & \cdots & d & -1 \\
            -1 & -1 & \cdots & -1 & d\\
            -1 & -1 & \cdots & -1 & -1 \\
        \end{bmatrix*},
\]
\[
    \GDn{d}:=
        \begin{bmatrix*}[r]
            -1 & 1 &  &  &  \\
            -1 & -1 & 1 &  & \\
                 &  & \ddots & \ddots &  &  \\
                 & & &  -1 & 1\\
                 & &  &    & -1 \\
        \end{bmatrix*},
        \text{ and }
    \GDndual{d}:=
        \begin{bmatrix*}[r]
            1 &  &    & & 1/2 \\
             & 1 &   & & 1/2 \\
              &  & \ddots &  & \multicolumn{1}{c}{\vdots} \\
              & &   & 1 & 1/2 \\
              & &   &  & 1/2 \\
        \end{bmatrix*}.
\]
\begin{table}[!htb]
\caption{Select \bxs s for
$d>3$. 
Shifts of the \bxs s for $\Z^d$, $\ZAn{d}$ and $\ZAndual{d}$ yield a basis.
$\DirSet{\Z^d}{2}$ are \bxs s on
$\ZDn{d}$.
        }
\begin{center}
\begin{tabular}{ccccccc}
\toprule
\multirow{2}{*}{lattice} & 
\multirow{2}{*}{dim.} & 
\multicolumn{2}{c}{direction sets } & 
\multirow{2}{*}{degree} & differentiability 
& \multirow{2}{*}{reference}
\\
&& 
$1$ & $2$  & & $r{-}2=$
\\
\midrule
\multirow{1}{*}{$\Z^d$} & $d{\ge}2$
& $n$ & $0$  & $d(n{-}1)$ &  $n{-}2$ &  B-splines \cite{deBoor1978Practical} 
\\
\cmidrule{1-7}

\multirow{1}{*}{$\ZAn{d}$} & $d{\ge}2$
& $1$
& 0 
& $d(d{-}1)/2$
&  $d{-}2$ 
&  \cite{Kim2011Symmetric}
\\
\cmidrule{1-7}

\multirow{1}{*}{$\ZAndual{d}$}& $d{\ge}2$
& $n$
& 0 
& $(d{+}1)n{-}d$
& $2(n{-}1)$ &  \cite{Kim2010Symmetric}
\\
\cmidrule{1-7}

\multirow{1}{*}{$\ZDn{d}$}& $d{\ge}2$
& $d(d{-}1)$ & 0  
& $d(d{-}2)$ 
& 
$
2d{-}4
$ 
&     \cite{Kim2011Symmetric}
\\
\cmidrule{1-7}

\multirow{4}{*}{$\ZDndual{d}$} & $d{=}4$
&  $1^\dag$ &  $0$ &  $4$ & $2$ & 
 \cite{Horacsek2022FastSpline}
\\
& $d{=}4$
& $1^{\dag\dag}$ & $0$ & $8$ & $4$ & \cite{Kim2011Symmetric}
\\
& $5{\le}d{\le}7$
& $1$ & $1$ & $2^{d-1}$ & $2^{d-2}$ & \cite{Kim2011Symmetric}
\\
& $d{>}4$
& $n$ & $0$ & $d(n-1)$ & $n-2$ & B-splines
\\
\bottomrule\multicolumn{6}{l}{\scriptsize 
$^\dag$  Constructed from
directions $\{((\pm1)^4)\}$ only. 
}
\\
\multicolumn{6}{l}{\scriptsize 
$^{\dag\dag}$  Constructed from directions
$\{\perm(\pm2,0^3)\}$ and $\{((\pm1)^4)\}$. }
\end{tabular}
\end{center}
        \label{tab:nvariate}
\end{table}

Note that $\GAn{d}$ and $\GAndual{d}$
are $(d+1)\times d$ and the corresponding lattices
are generated in the hyperplane 
of the equation $x_1+\cdots+x_{d+1}=0$.

\tabref{tab:DS.nvariate} lists the first and second
direction sets of the five lattices.
As in the bi- and the trivariate cases,
various symmetric \bxs s can be constructed from
these directions.
We observe that for $\ZDndual{d}$ , $d>4$, there is a rich set of first
directions, 
all corresponding to B-splines,
to build smooth symmetric splines.
\tabref{tab:nvariate} lists 
%
some important classes of \bxs s whose shifts live on these high-dimensional lattices, see e.g.\ \cite{Kim2011Symmetric}.
Note that for some dimensions,
two different direction sets share the same distance:
for $\ZDndual{4}$ there are $16/2+8/2=12$ first directions 
of the patterns
$(\pm1,\pm1,\ldots,\pm1)$
and $\pi(\pm2,0,0,0)$
and either or both
groups yields a symmetric \bxs.

\begin{figure}[ht]
    \includegraphics[width=.2\textwidth]{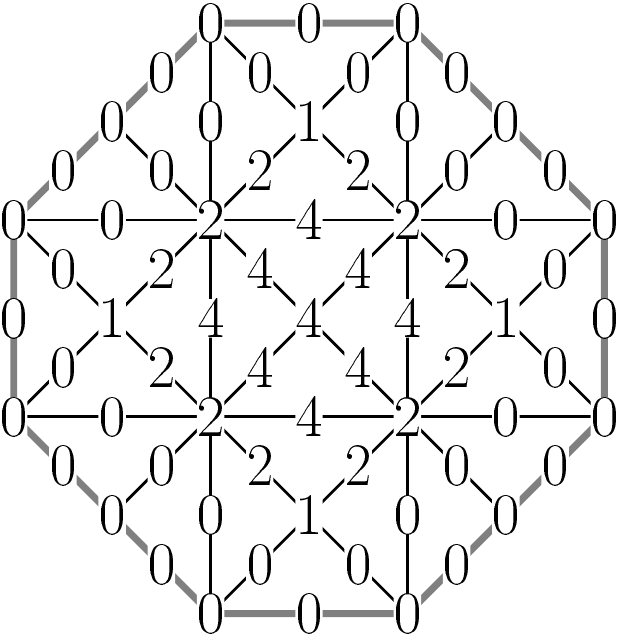}
    \caption{The polynomial pieces in the support of $\Mcc{11}$ and
    the BB-net (scaled by $8$).}
    \label{fig:polyc11}
\end{figure}

\begin{figure}[ht]
	\begin{tabular}{cc}
		\subfigure[]{
			\label{fig:Mc21.pieces}
		\adjustimage{width=.2\textwidth,valign=m}{Figures/Z2/Mc21-pieces}	
		}
&
	\begin{tabular}{cccc}
		\subfigure[]{%
			\label{fig:Mc21.bb.1}%
		\includegraphics[width=.15\textwidth]{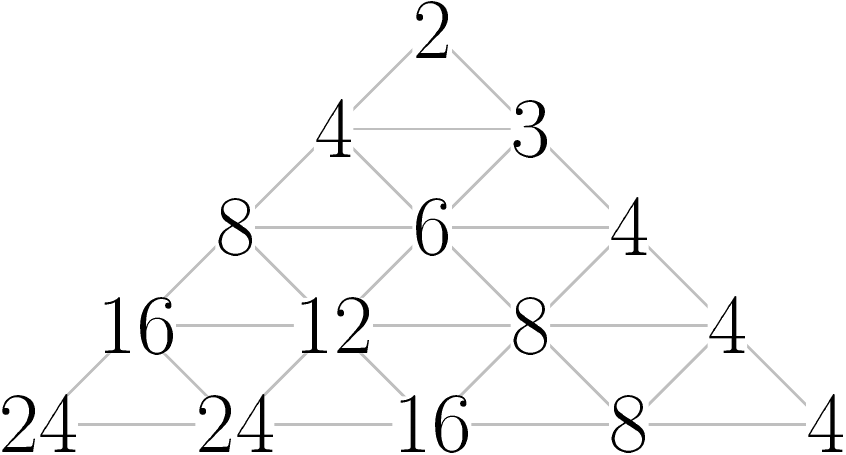}	%
		}	&
		\subfigure[]{%
			\label{fig:Mc21.bb.2}%
		\includegraphics[width=.15\textwidth]{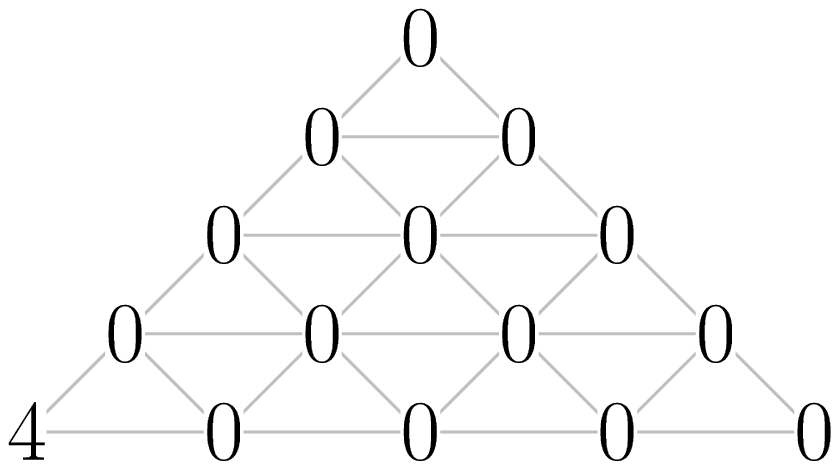}	%
		}	&
		\subfigure[]{%
			\label{fig:Mc21.bb.3}%
		\includegraphics[width=.15\textwidth]{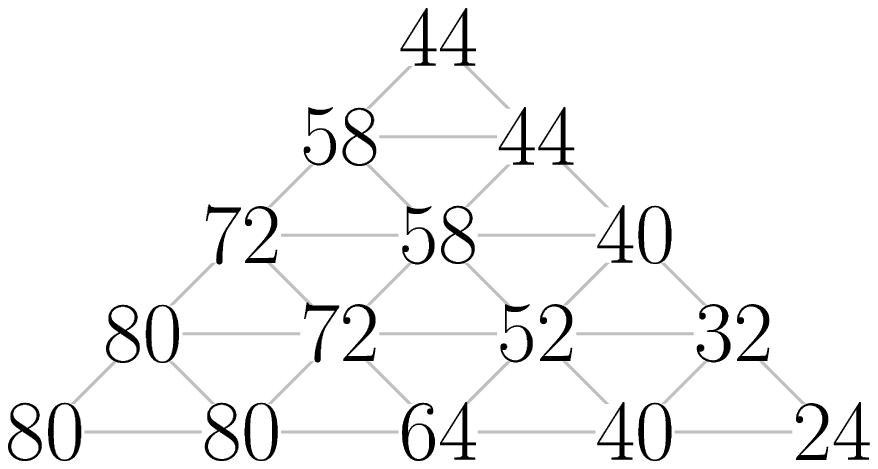}	%
		}	&
		\subfigure[]{%
			\label{fig:Mc21.bb.4}%
		\includegraphics[width=.15\textwidth]{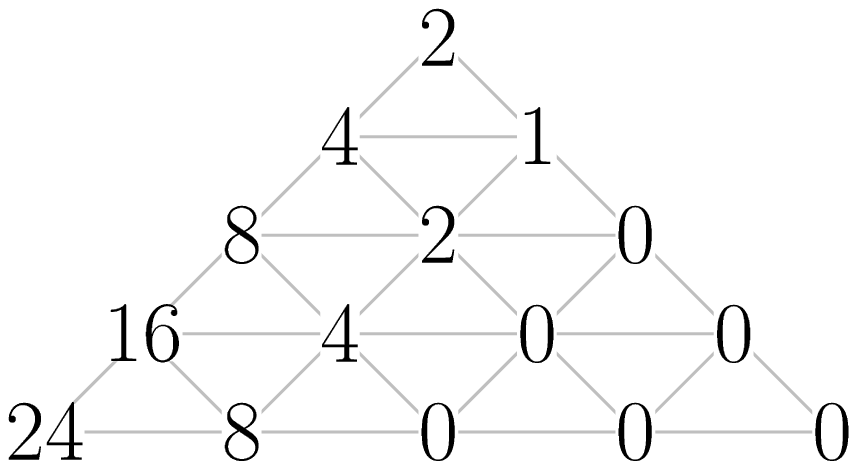}	%
		}	\\
		\subfigure[]{%
			\label{fig:Mc21.bb.5}%
		\includegraphics[width=.15\textwidth]{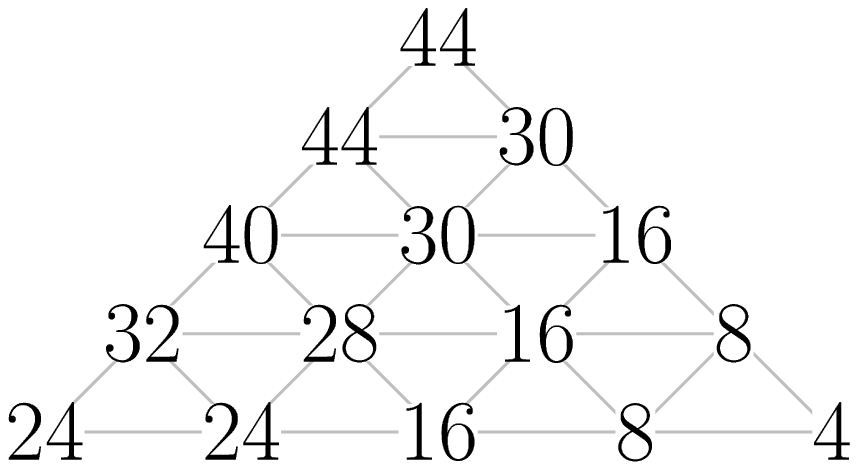}	%
		}	&
		\subfigure[]{%
			\label{fig:Mc21.bb.6}%
		\includegraphics[width=.15\textwidth]{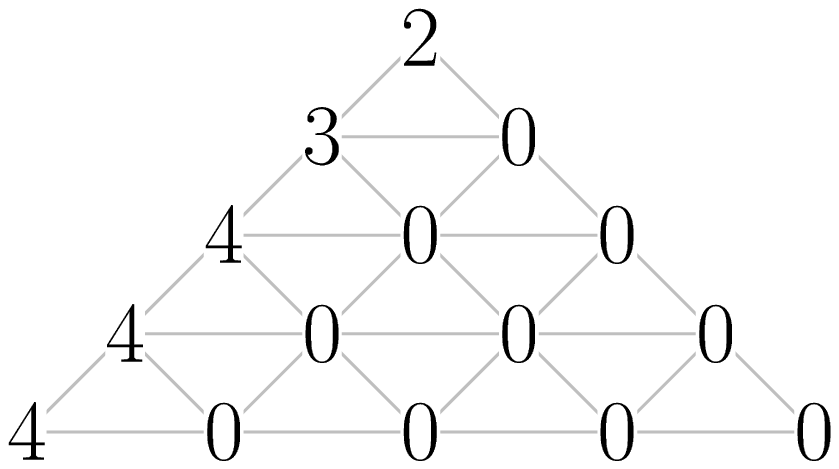}	%
		}	&
		\subfigure[]{%
			\label{fig:Mc21.bb.7}%
		\includegraphics[width=.15\textwidth]{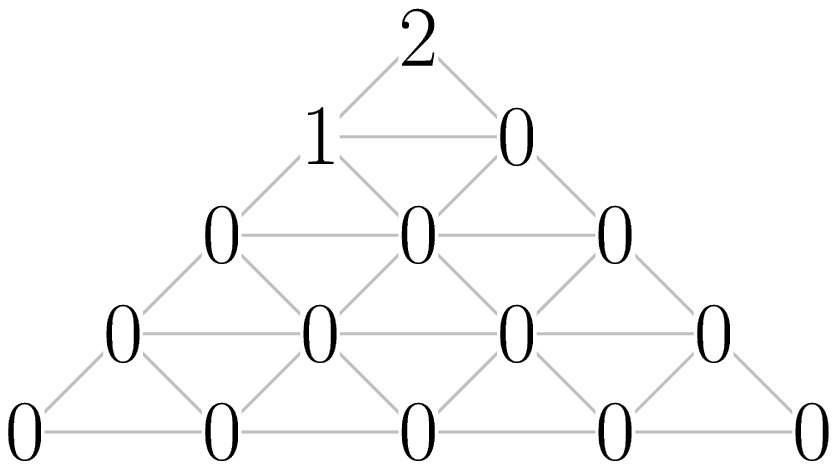}	%
		}%
	\end{tabular}
\end{tabular}
\caption[]{
    From \cite{Chui1991Algorithms}.
\subref{fig:Mc21.pieces} The polynomial pieces in the support of $\Mxzp$.
Pieces of the same color
have the same BB-net after appropriate rigid transformation and the 
BB-nets  (multiplied by $192$) of the pieces labeled b,$\ldots$,h are shown in
\subref{fig:Mc21.bb.1}--\subref{fig:Mc21.bb.7}.
	
}
\label{fig:polyc21}
\end{figure}
\begin{figure}[ht]
    \centering
    \begin{tabular}{cccc}
    \subfigure[]{
    \label{fig:Mh20.pieces}
    \includegraphics[width=.2\textwidth]{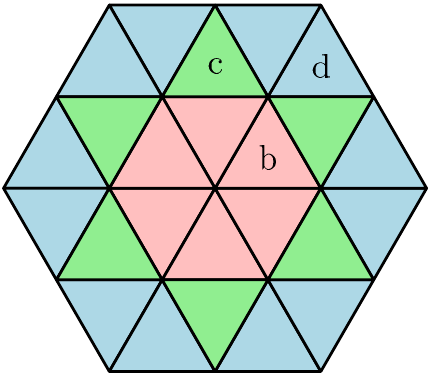}
    }
    &
    \subfigure[]{
    \label{fig:Mh20.bb.1}
    \includegraphics[width=.2\textwidth]{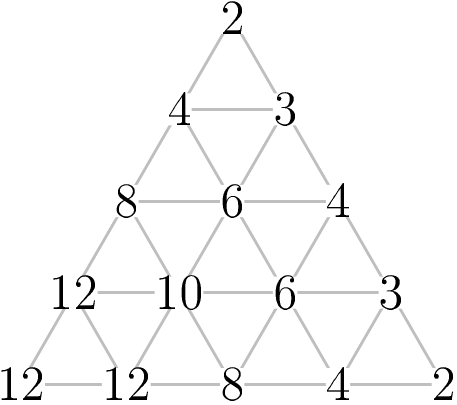}
    }
    &
    \subfigure[]{
    \label{fig:Mh20.bb.2}
    \includegraphics[width=.2\textwidth]{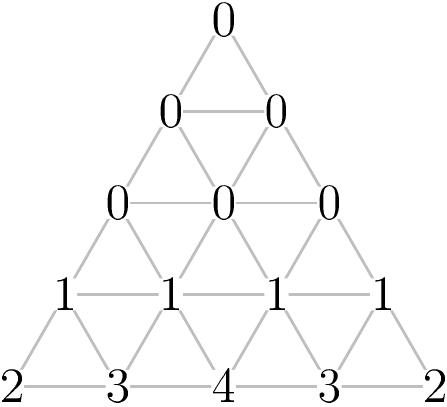}
    }
    &
    \subfigure[]{
    \label{fig:Mh20.bb.3}
    \includegraphics[width=.2\textwidth]{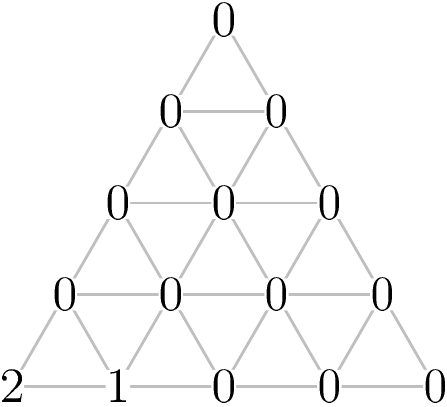}
    }
    \end{tabular}
    \caption{
    From \cite{Chui1991Algorithms}.
\subref{fig:Mh20.pieces} The polynomial pieces in the support of
$\MhexII$. 
Pieces of the same color
have the same BB-net after appropriate rigid transformation and the 
BB-nets  (multiplied by $24$) of the pieces labeled b,c,d are shown in	\subref{fig:Mh20.bb.1},\subref{fig:Mh20.bb.2},\subref{fig:Mh20.bb.3}.
}
    \label{fig:Mh20.bb}
    \label{fig:polyh20}
\end{figure}

\section{Conversion to piecewise polynomial form}

It is useful to express the \bxs\ pieces as polynomials, and
in particular in the Bernstein-B\'ezier (BB-) form, see e.g.\
\cite{deBoor1987B-form}.
The partition into pieces follows from the 
convolution  directions.  
The BB-coefficients  are obtained 
from the differentiability constraints across boundaries
and by normalizing the map, see \cite{Lai2007Spline}.
\figref{fig:polyc11},
\figref{fig:polyc21} and \figref{fig:polyh20}
show examples of the re-representation in BB-form.
For trivariate \bxs s, using the 
the constraints can be error-prone.
An easier approach is to sample the spline
at sufficiently many interior points, 
using one of \cite{deBoor1993Evaluation,Kobbelt1997Stable},
and solve
for the BB-coefficients, keeping in mind that 
the coefficients are integers after scaling by a known multiple (see  \cite{Kim2009Fast}); or, and this is  faster and yields polynomial pieces in partially factored form,
to apply a Green's function decomposition
and inverse Fourier transform
\cite{Horacsek2018Closed}.



\section{Efficient evaluation}
\label{sec:eval}
By reversing the convolution, the algorithms of
\cite{deBoor1993Evaluation,Kobbelt1997Stable}
evaluate \bxs s recursively. This process is stable
except near the boundaries between the polynomial pieces (knot lines in 2D, knot planes in 3D).
Near boundaries \cite{deBoor1993Evaluation} applies random perturbation and
\cite{Kobbelt1997Stable} propose careful bookkeeping.
Converting
the \bxs\ pieces to  BB-form
yields
much faster and stable evaluation \cite{Kim2009Fast}, also of derviatives.
A general technique to accelerate evaluation
is to leverage symmetry 
\cite{Kim2017Analysis,Horacsek2022FastSpline} with a general implementation available at 
\cite{Horacsek2022FastSpline} that automates steps and 
    generates
    GPU kernels.
\tabref{tab:eval} lists \bxs s with an available
optimized evaluation code, 
some implemented on the GPU for high parallelism.
\begin{table}[!htb]
    \caption{Some fast  3D \bxs\  evaluation implementations.
    See also \cite{Horacsek2022FastSpline}.
    }
    \label{tab:eval}
    \centering
    \begin{tabular}{ccc}
    \toprule
    \bxs\ & algorithm & code \\
    \midrule
    $\Mcc{400}$
    & \cite{Sigg2005Fast,Kim2017Efficient} 
    & \cite{Sigg2005Fast}
    \\
     $\Mcc{010}$ &
     \cite{Kim2023Volume} &
      \cite{Kim2023Volume}
       \\
     $\Mcc{101}$ & \cite{Kim2021Fast} & 
     \cite{Kim2021Fast}
     \\
$\Mfcc{100}$ & \cite{Kim2008Box,Kim2013GPU} &
 \cite{Kim2013GPU} \\
 $\Mbcc{200}$ & \cite{Casciola2006Mixed,Finkbeiner2010Efficient} 
 & \cite{Finkbeiner2010Efficient} \\
 $\Mbcc{110}$ & \cite{Kim2013Quartic,Kim2012BCC,Kim2017Analysis} & 
\cite{bcc7}
\\
$\Mbcc{040}$ & \cite{Csebfalvi2006Prefiltered} 
& \\ \bottomrule
    \end{tabular}

\end{table}
Subdivision offers a stable and fast alternative when 
rendering an approximation,
say a triangulation 
of a bivariate \bxs\ graph. An alternative
approximate evaluation is based on Fast Fourier transform \cite{McCool1996Accelerated}.

\section{Use for reconstruction or approximation}
\label{sec:reconstruction}



\begin{figure*}[ht]
\begin{center}
\def\scalewidth{.24}
    \subfigure[$\Mcc{300}$]{%
    \includegraphics[width=\scalewidth\textwidth]{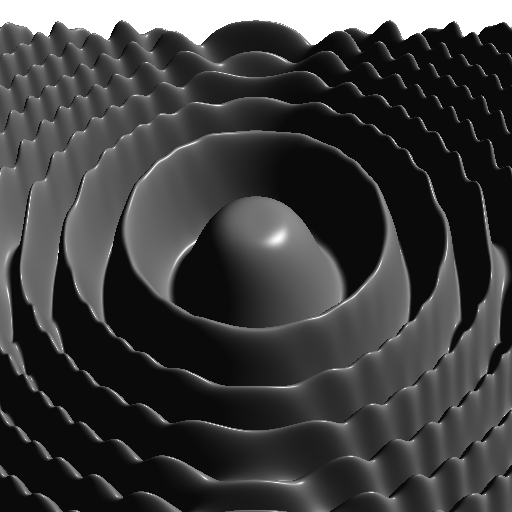}%
        }%
    \subfigure[$\Mcc{010}$]{%
    \includegraphics[width=\scalewidth\textwidth]{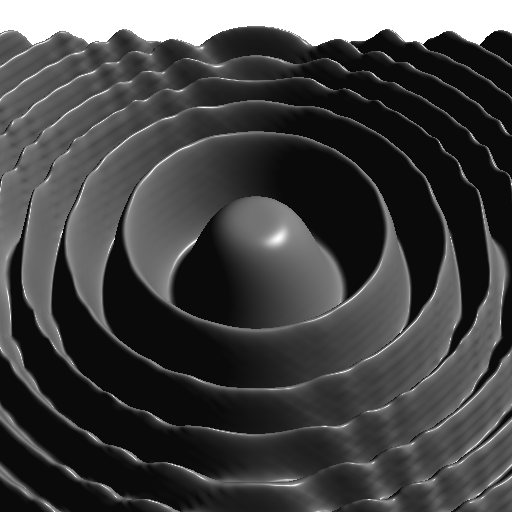}%
        }%
    \subfigure[$\Mcc{400}$]{%
    \includegraphics[width=\scalewidth\textwidth]{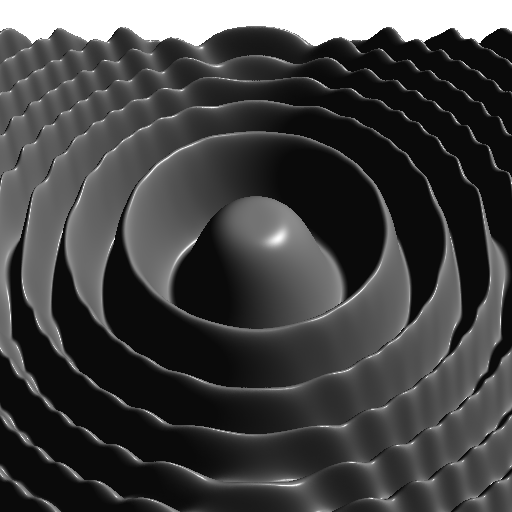}%
        }%
    \subfigure[$\Mcc{101}$]{%
    \includegraphics[width=\scalewidth\textwidth]{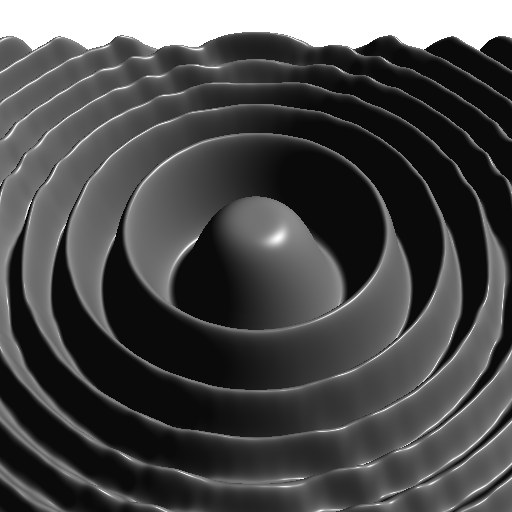}%
        }
    \subfigure[$\Mfcc{100}$]{%
    \includegraphics[width=\scalewidth\textwidth]{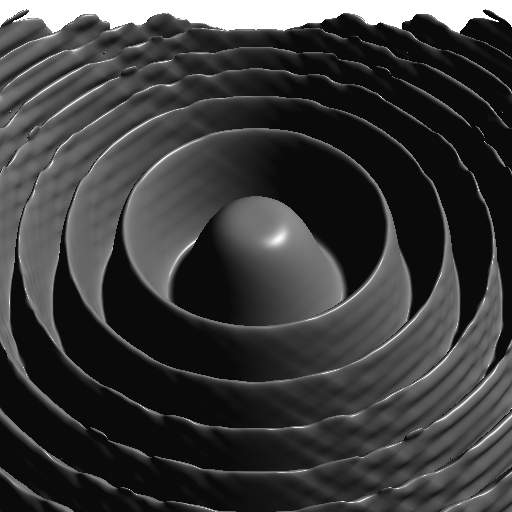}%
        }%
    \subfigure[$\Mbcc{200}$]{%
    \includegraphics[width=\scalewidth\textwidth]{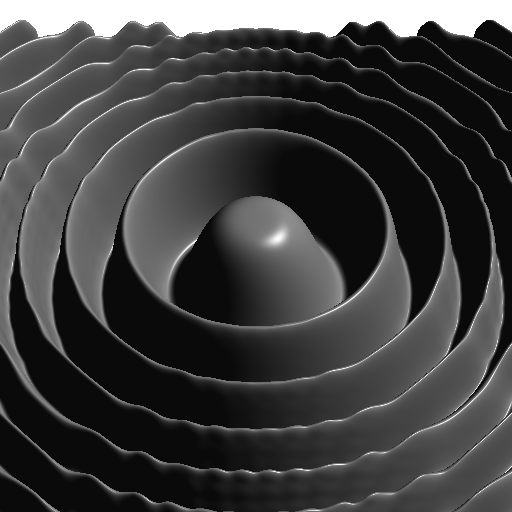}%
        }%
    \subfigure[$\Mbcc{110}$]{%
    \includegraphics[width=\scalewidth\textwidth]{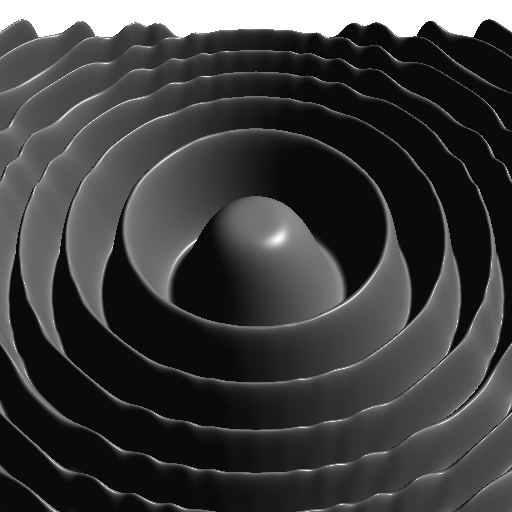}%
        }%
    \subfigure[ground truth]{%
    \label{fig:ML.truth}%
    \includegraphics[width=\scalewidth\textwidth]{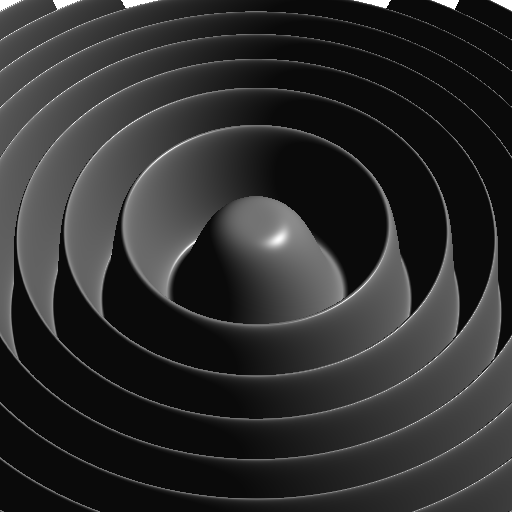}%
        }%
%
\end{center}
\caption{
Ray-intersection rendering (ray-casting) of a level set of the Marschner-Lobb signal \subref{fig:ML.truth}
with identical
sampling density on their domain lattices.
}
\label{fig:ML}
\end{figure*}

\begin{table}[!htb]
    \caption{Quasi-interpolants 
    of select \bxs s
   of approximation order (a.o.) $3$ or $4$.
   Note that $q_0$ and $q_1$ are scaled
   for clearer presentation.
    }
    \label{tab:qi}
    \centering
    \begin{tabular}{cccccc}
        \toprule
        lattice
        & a.o.
        & \bxs         
        & $24q_0$ 
        & $-12q_1$ 
        & references 
        \\
        \midrule
        \multirow{2}{*}{$\Z^2$}
        & $3$
        &
        $\Mcc{30}$, $\Mcc{11}$ 
        &         $18$ & $3$
        & \cite{Kim2011Symmetric}
        \\
        \cmidrule{2-6}
        & $4$ 
        & $\Mcc{20}$, $\Mcc{21}$ 
        & $40$ & $4$
        \\
        \cmidrule{1-6}
        $\Zhex$ 
        & $4$
        & $\Mhex{20}$ 
        & $13$ & $2$
        & \cite{Kim2010Symmetric}
        \\
        \cmidrule{1-6}
        \multirow{5}{*}{$\Z^3$} 
        & \multirow{2}{*}{$3$}
        & $\Mcc{300}$  
        & $21$ & $3$
        \\
        && $\Mcc{010}$  
        & $24$ & $4$
        & \cite{Kim2023Volume}
        \\
        \cmidrule{2-6}
        & \multirow{3}{*}{$4$}
        & $\Mcc{400}$  
        & $24$ & $4$
        \\
        & & $\Mcc{101}$ 
        & $27$ & $5$
        & \cite{Kim2011Symmetric}
        \\
        & & $\Mcc{002}$ 
        & $36$ & $8$
        \\
        \cmidrule{1-6}
        \multirow{3}{*}{$\Zfcc$} 
        & \multirow{2}{*}{$3$}
        &
        $\Mfcc{100}$  
        & $18$ & $1$
        & \cite{Kim2011Symmetric}
            \\
        &&     $\Mfcc{030}$  
        & $30$ & $3$
            \\
            \cmidrule{2-6}
            & $4$ 
        &     $\Mfcc{040}$ 
        & $36$ & $4$
        \\
        \cmidrule{1-6}
        \multirow{4}{*}{$\Zbcc$} 
        & \multirow{2}{*}{$3$}
        &
        $\Mbcc{030}$ 
        & $24$ & $3$
            \\
        & & $\Mbcc{001}$  
        & $28$ & $4$
            \\
            \cmidrule{2-6}
            & \multirow{2}{*}{$4$}
        & $\Mbcc{200}$, $\Mbcc{110}$ 
        & $20$ & $2$
            & \cite{Entezari2009Quasi,Kim2010Symmetric,Kim2013Quartic}
            \\
        & & $\Mbcc{040}$ 
        & $28$ & $4$
        \\
        \bottomrule
    \end{tabular}

\end{table}
A promising application of \bxs s is the
approximation and reconstruction of a function $f$
from samples $\{f(\vj):\vj\in\ZG\}$ on a lattice
$\ZG$.
To attain the maximal approximation order
of the \bxs\ space,
i.e., to obtain $c$ in Eq. \eqref{eq:a.o},
the samples
are 
convolved with a discrete \emph{quasi-interpolant} to form the control points
\[
      c(\vj):= 
        q_0
      f(\vj)
      + 
      q_1
      \sum_{\vk\in\DirSet{\ZG}{1}}
      \left(f(\vj+\vk) 
        + f(\vj-\vk)\right),
        \quad\forall \vj\in\ZG
\]
of the optimally approximating spline 
$
     \sum_{\vj\in\ZG}c(\vj)M(\cdot-\vj).
$
Several techniques exist to derive 
 quasi-interpolants for \bxs s. 
\cite{deBoor1973Spline,deBoor1993Box,Blu1999Quantitative-I,Condat2007Quasi}.
\tabref{tab:qi} lists  quasi-interpolants,
defined by $q_0$ and $q_1$,
for the \bxs s of approximation order $3$ or $4$
of \tabref{tab:bivariate} and \tabref{tab:trivariate}.
Level sets of quasi-interpolating functions in three variables
are used to display
Computed Tomography (CT)
and Magnetic Resonance Imaging (MRI)
data. 
A standard test function 
is the 
 Marschner-Lobb  signal
\cite{Marschner1994Evaluation},
a combination of Dirac pulses and a circularly symmetric, disc-shaped component, see \figref{fig:ML.truth}.
\figref{fig:ML} compares how 
convolution directions enhance or prevent reproduction of the circular features.

\section{Splines from pieces and unions of boxes}
One can consider the characteristic 
function of a piece of the  box
or of a union of boxes,
and then convolve these characteristic functions.
%
Convolving the characteristic function of half of a box in 2D, i.e.\
of a triangle, 
yields 
\emph{half-box spline} spaces with properties akin to \bxs s 
\cite{Sabin1977Use,Farin1982Desining,Boehm1987On,Prautzsch2002Box,Barendrecht2019Bivariate}.
%
%
%
%
Alternatively, one can
juxtapose non-centered boxes to form the  Voronoi  cell of a lattice, i.e.\ the region nearest to a lattice point.
The convolution of the characteristic function of the Voronoi cell then yields
\emph{Voronoi splines}
\cite{VanDeVille2004Hex,mirzargar2010voronoi}. 
Voronoi splines provide an example of how asymmetric splines can be linearly combined
to form symmetric splines.
Note though that such splines 
typically do not yield nested spaces
\cite{Peters2014Refinability}.



\section{Conclusion}
\label{sec:con}
Symmetric \bxs s provide a mature and powerful framework for  shift-invariant  smooth functions  on a lattice. 
For bi- and tri-variate splines, a number of 
efficient \bxs s are now well-documented
and come with optimized implementations. 

\bigskip
\noindent
\textbf{Acknowledgements}
This work was supported by a 2022 sabbatical year research grant of the University of Seoul hosted by the University of Florida.
We thank Carl de Boor for feedback on an early draft.
\bibliographystyle{abbrvnat}
\bibliography{p.bib}

\end{document}